\documentclass[12pt]{amsart}
\usepackage{graphicx}
\usepackage{amssymb}
\usepackage{url}
\usepackage{amsmath}
\usepackage{latexsym}
\usepackage{verbatim,euscript}
\usepackage{fullpage,color}

\usepackage[matrix,arrow,curve]{xy}
\usepackage[small,nohug,heads=vee]{diagrams}
\diagramstyle[lablestyle=\scriptstyle]
\usepackage{lmodern}
\DeclareMathSizes{7}{7}{5}{4} % \scriptsize is 7 pt
\DeclareMathSizes{5}{5}{3}{2} % \tiny is 5 pt

\newcommand{\Lg}{\mbox{$\mathfrak g$}}

\newcommand{\Lh}{\mbox{$\mathfrak h$}}
\newcommand{\Lk}{\mbox{$\mathfrak k$}}
\newcommand{\Lp}{\mbox{$\mathfrak p$}}
\newcommand{\La}{\mbox{$\mathfrak a$}}

\newcommand{\Lt}{\mbox{$\mathfrak t$}}
\newcommand{\Lq}{\mbox{$\mathfrak q$}}
\newcommand{\Ls}{\mbox{$\mathfrak s$}}

\newcommand{\Ll}{\mbox{$\mathfrak l$}}

\newcommand{\g}[1]{\mbox{$\mathfrak{#1}$}}

\newcommand{\Pf}{{\em Proof}. }
\newcommand{\EPf}{\hfill$\square$}

\newcommand{\C}{\mbox{$\mathbb C$}}

\newcommand{\cat}{\mbox{$/\!\!/ $}}

\usepackage[colorlinks=true,linkcolor=blue,urlcolor=my_color,citecolor=magenta]{hyperref}

\definecolor{my_color}{rgb}{0,0.5,0.5}
\definecolor{MIXT}{rgb}{0.4,0.3,0.6}

\newtheorem{thm}{Theorem}

\newtheorem{cor}[thm]{Corollary}
\newtheorem{prop}[thm]{Proposition}
\newtheorem{lem}[thm]{Lemma}

\theoremstyle{remark}
\newtheorem{rmk}[thm]{Remark}

\theoremstyle{definition}
\newtheorem{ex}[thm]{Example}

\numberwithin{equation}{section}

\begin{document}

\title{Polar symplectic representations}

\author{Laura Geatti}          

\address{
         Dipartimento di Matematica, Universit\`a di Roma 2 Tor Vergata,
         via della Ricerca Scientifica, 00133 Roma, Italy}
\email{geatti@mat.uniroma2.it}

\author{Claudio Gorodski}\thanks{The second author has been partially supported by the CNPq grant 303038/2013-6 and the FAPESP project 2011/21362-2.}

\address{     Instituto de Matem\' atica e Estat\'\i stica,
         Universidade de S\~ ao Paulo,
         Rua do Mat\~ao, 1010,
         S\~ ao Paulo, SP 05508-090,
         Brazil}
\email{gorodski@ime.usp.br}



\maketitle

\begin{abstract}
We study polar representations in the sense of Dadok and Kac 
which are symplectic. We show that such representations are coisotropic
and use this fact to give a classification. We also study their moment maps
and prove that they separate closed orbits. 
Our work can also be seen as a specialization of some 
of the results of Knop on multiplicity free symplectic 
representations to the polar case.
\end{abstract}

\section{Introduction}

A rational representation of a complex reductive linear algebraic
group~$G$ on a finite-dimensional complex vector space~$V$ is called
\emph{polar} if there exists a subspace $c\subset V$ consisting
of semisimple elements such that $\dim c = \dim V\cat G$ (the categorical
quotient), and for a dense subset of $c$, 
the tangent spaces to the orbits are 
parallel~\cite{DK}; then it turns out that every closed orbit 
of $G$ meets $c$ (Prop.~2.2, \emph{ibid}).
In this paper we study the class of polar representations 
which are \emph{symplectic}, namely, preserve a non-degenerate
skew-symmetric bilinear form~$\omega$ on~$V$ (polarity of a 
representation depends only on the identity component, and we 
assume throughout that all groups are connected).  
We first prove:

\begin{thm}\label{polar-implies-coisotropic}
A polar symplectic representation is coisotropic.
\end{thm}

Recall that a symplectic representation $V$ of $G$ is 
\emph{coisotropic} if a generic $G$-orbit is co\-isotropic, namely, 
$(\mathfrak g\cdot v)^{\omega} \subset \mathfrak g\cdot v$ where 
$v\in V$ is generic, $\mathfrak g$ denotes the Lie algebra of~$G$
and $()^\omega$ refers to the symplectic complement.  
Representations in this class can be characterized by a 
number of different properties, e.g.~the Poisson algebra
of invariants $\mathbb C[V]^G$ is commutative
(cf.~\cite[p.~224 and Prop.~9.1]{Knop2} and~\cite[Introd.]{Lo});
in particular, they are also called \emph{multiplicity-free}
(in the symplectic sense).

Using Theorem~\ref{polar-implies-coisotropic}, we can
reduce the classification of polar symplectic representations,
up to geometric equivalence, 
to that of coisotropic representations given in~\cite{Knop1}.
In contrast to the case of coisotropic representations, 
it turns out that every saturated decomposable polar symplectic
representation is an outer product (see section~\ref{prelim} for 
unexplained terminology).

\begin{thm}\label{classif}
The saturated indecomposable polar symplectic representations are listed
in Tables~A and~B. 
Every saturated polar symplectic 
representation is an outer product of 
indecomposable polar symplectic representations.
\end{thm}

\[\begin{array}{|c|c|c|c|}
\hline
G & V & \textsl{$\dim V\cat G$} & \textsl{Conditions} \\
\hline
SO_p\otimes Sp_{2m}& \C^p\otimes \C^{2m} & \min\{[\frac p2],m\} & m\geq1,\ p\geq3\\
Sp_{2m} &\C^{2m} & 0 & m\geq1 \\
SL_2\times Spin_7 &\C^2\otimes \C^8  & 1 & - \\
SL_2\times Spin_9  &\C^2\otimes \C^{16} & 2 & - \\
Spin_{11} &\C^{32} & 1 & - \\
Spin_{12} &\C^{32} & 1 & - \\
Spin_{13} &\C^{64} & 2 & - \\
SL_2 & S^3(\C^2) & 1 & - \\
SL_6 &\Lambda^3(\C^6) & 1 & - \\
Sp_6 &\Lambda^3(\C^6)\ominus\C^6 & 1 & - \\
SL_2\times G_2 & \C^2\otimes \C^7 & 1 & - \\
E_7 & \C^{56} & 1  & - \\
\hline
\end{array} \]
\begin{center}
\textsc{Table A: Indecomposable polar symplectic representations
of type~1 }
\end{center}

\[\begin{array}{|c|c|c|c|c|}
\hline
G & V & \textsl{$\dim V\cat G$} & \textsl{Conditions} & \textsl{$\C^\times$ essential}\\
\hline
\C^\times \times SL_m\times SL_n& \C^m\otimes \C^n\oplus(\C^m\otimes\C^n)^* 
& n & m\geq n\geq2 & \textrm{yes iff $m=n$}\\
GL_n &\Lambda^2(\C^n)\oplus\Lambda^2(\C^n)^* & [\frac n2] & n\geq4 & \textrm{yes iff $n$ even}\\
GL_n &S^2(\C^n)\oplus S^2(\C^n)^* & n & n\geq2 & \textrm{yes}\\
GL_n &\C^n\oplus\C^{n*} & 1 & n\geq1 & \textrm{yes iff $n=1$}\\
\C^\times\times Sp_{2m} &\C^{2m}\oplus\C^{2m*}&  1 & m\geq2 & \textrm{no}\\
\C^\times\times SO_m & \C^m\oplus\C^{m*} & 2 & m\geq5 &\textrm{yes}\\
\C^\times\times Spin_7 &\C^8\oplus\C^{8*} & 2 & - &\textrm{yes}\\
\C^\times\times Spin_{10} &\C^{16}\oplus\C^{16*} & 2 & - &\textrm{no}\\
\C^\times\times G_2 & \C^7\oplus\C^{7*} & 2 & - & \textrm{yes}\\
\C^\times\times E_6 & \C^{27}\oplus\C^{27*} & 3  & - & \textrm{yes}\\
\hline
\end{array} \]
\begin{center}
\textsc{Table B: Indecomposable polar symplectic representations
of type~2}
\end{center}

In the last column of Table~B, non-essentialness of the center means that 
its removal does not change the closed orbits; otherwise, the closed
orbits change and the representation ceases to be polar. 

A \emph{symplectic symmetric space} is a symmetric space which is endowed
with a symplectic structure invariant by the symmetries. 
Our interest in them is that the (complexified) isotropy representations of 
symplectic symmetric spaces provide examples of symplectic 
$\theta$-groups~\cite{V,Kac2},
thus, polar symplectic representations. Conversely, it is a natural question 
to ask which polar symplectic representations come from 
symplectic symmetric spaces.
We say that two symplectic representations are 
\emph{closed orbit equivalent} if there exists a symplectic isomorphism 
between the representation spaces mapping closed orbits onto closed orbits
(for the sake of comparison, recall that in the orthogonal case
all polar representations come from symmetric spaces, up to closed
orbit equivalence~\cite{D,GG}). Note that polarity is a property 
of closed orbit equivalence classes. 

\begin{thm}\label{symmetric-spaces}
A polar symplectic representation is closed orbit equivalent 
to the isotropy representation of a complex semisimple 
symplectic symmetric space
if and only if it is closed orbit equivalent to the complexification of the 
isotropy representation of a semisimple Hermitian Riemannian symmetric space.
In the saturated case, such representations are exactly the outer products
of representations listed in Table~B.  
\end{thm} 

It is relevant to notice that the only cases in Table~A 
which are not $\theta$-groups are the representations of
$SL_2\times Spin_9$, $Spin_{11}$, $Spin_{13}$ and $SL_2\times G_2$~\cite{DK,L},
and that only two of them have $\dim V\cat G>1$.

Finally, recall that a symplectic representation $(G,V)$
has a canonical \emph{moment map} $\mu:V\to\Lg^*$ 
(see section~\ref{sec:moment}). 
Since $\mu$ is equivariant, it induces an \emph{invariant moment
map} $\psi=\mu\cat G:V\cat G\to\Lg^*\cat G$.

\begin{thm}\label{moment-map}
The moment map of a 
saturated polar symplectic representation
maps closed orbits to closed orbits and separates closed orbits.
\end{thm}

\begin{rmk}
The only saturated indecomposable polar symplectic representation
for which the invariant moment map $\psi$ fails to be
an isomorphism from
$V\cat G$ onto an affine space in $\Lg^*\cat G$ is the last 
one in Table B.
Hence, in all the other cases the morphism 
$\psi^*:\mathbb C[\Lg^*]^G\to\mathbb C[V]^G$ is surjective,
that is, all invariants are pull-backs of coadjoint invariants.
\end{rmk}
 
\begin{rmk}
In case of type~2 representations, Theorems~\ref{polar-implies-coisotropic}
and~\ref{moment-map} reduce to known facts about polar representations
of compact Lie groups in the sense of Dadok~\cite{D}. 
Let $(K,U)$ be an orthogonal representation of a compact Lie group $K$
and consider its complexification $(G:=K^{\mathbb C},V:=U^{\mathbb C})$. 
It is easy to check that $(G,V)$ is polar if and only if
$(K,U)$ is polar. Suppose now $U$ admits an invariant complex structure. 
Then $(G,V=U\oplus U^*)$ is coisotropic if and only 
$(K,U)$ is multiplicity free (cf.~\cite[p.532]{Knop1} 
or~\cite[Prop.~9.2]{Knop2})
if and only if $(K,U)$ has coisotropic principal $K$-orbits
(\cite[Thm.~3.1]{Knop3} and \cite[Prop.~12]{V2})
%(\cite[Thm.]{DP} and \cite[p.274]{HW}), 
so Theorem~\ref{polar-implies-coisotropic} says that a polar 
representation of a compact Lie group has coisotropic principal
$K$-orbits (compare~\cite[Thm.~1.1 and Lem.~2.7]{PTh2}). Moreover the moment
map $\mu$ of $(G,V)$ restricts to the moment map $\mu_K$ of $(K,U)$,
every $G$-orbit through a point in $U$ is closed~\cite{Bi}, and
two different $K$-orbits in~$U$ cannot be contained in the same
$G$-orbit (since $\mathbb C[V]^G=\mathbb R[U]^K\otimes\mathbb C$;
see also~\cite[\S~2.2]{Br}), so Theorem~\ref{moment-map} says that  
$\mu_K$ separates $K$-orbits (compare~\cite[Cor.~1.5]{PTh2} and~\cite[p.~274]{HW}). 
\end{rmk}

The authors wish to thank Friedrich Knop, Paul Levy and the anonymous
referee for their valuable comments and suggestions which have 
substantially helped improve this work.

\section{Preliminaries}\label{prelim}

We begin by recalling terminology from~\cite{Knop2} that will be useful
in the sequel.
A symplectic representation of $G$ is called \emph{indecomposable} if it is not
isomorphic to the sum of two non-trivial symplectic representations of $G$. 
A symplectic representation $V$ of $G$ is called of \emph{type~1}
if $V$ is irreducible as a $G$-module, and it is called 
of~\emph{type~2} if $V=U\oplus U^*$ where $U$ is an irreducible
$G$-module not admitting a symplectic structure and the symplectic 
form on $V$ is given (up to a multiple) by 
\[ \omega(u_1+u_1^*,u_2+u_2^*)=u_1^*(u_2)-u_2^*(u_1). \]
Every indecomposable symplectic representation is either of type~1 or~2.
Two symplectic representations are isomorphic as $G$-modules if and only 
if they are isomorphic as symplectic representations. Every symplectic 
representation is a direct sum of finitely many indecomposable symplectic 
representations, and the summands are unique up to 
permutation~\cite[Thm.~2.1]{Knop1}.

It is convenient to revisit the result above as follows. 
Choose a maximal compact subgroup $K$ of $G$ (necessarily connected)
and a $K$-invariant Hermitian inner product $h$ on $V$. A $K$-invariant
conjugate linear automorphism $\epsilon:V\to V$ is then defined by
\begin{equation}\label{omega-h}
 \omega(u,v) = h (u,\epsilon v) 
\end{equation}
for all $u$, $v\in V$. Then 
\begin{eqnarray*}
h(u,\epsilon^2 v) &=&  \omega(u,\epsilon v) 
= -\omega(\epsilon v,u) 
= -h (\epsilon v,\epsilon u) \\
&=& -\overline{h(\epsilon u,\epsilon v)} 
= -\overline{\omega(\epsilon u,v)} \
= \overline{\omega(v,\epsilon u)} 
=\overline{h(v,\epsilon^2 u )} 
= h(\epsilon^2 u,v), 
\end{eqnarray*}
so $\epsilon^2$ is a $\C$-linear $K$-invariant Hermitian endomorphism of $V$.
It also follows from the above that $h(u,\epsilon^2 u )=-||\epsilon u||^2$, 
so $\epsilon^2$ is negative definite. Now there is a $h$-orthogonal 
$K$-irreducible decomposition $V=\bigoplus V_j$ such that 
$\epsilon^2|_{V_j}=\lambda_j\mathrm{id}_{V_j}$ for $\lambda_j<0$ and all~$j$. 
For each $j$, either $\omega|_{V_j\times V_j}$ is non-degenerate or it is zero
(since $\omega$ is $K$-invariant and $V_j$ is $K$-irreducible). 
In the former case, $\epsilon(V_j)=V_j$. In the latter case, 
$\epsilon(V_j)\perp_h V_j$ and $\epsilon(V_j)=V_j^*$ (since 
$\epsilon$ is conjugate-linear). Hence 
$V$ is an $h$-orthogonal direct sum of symplectic 
representations of \emph{type~1} ($V_j$ is irreducible and anisotropic)
and \emph{type~2} ($V_j\oplus V_j^*$, where $V_j$ is irreducible and 
isotropic). By renormalizing $h$, we may assume that 
$\epsilon^2=-\mathrm{id}_V$; in particular, $\epsilon$ becomes
an $h$-isometry.

Let $\rho_i:\Lg_i\to\mathfrak{sp}(V_i)$ for $i=1$, $2$ be two 
symplectic representations. We say $V_1$ and $V_2$ are 
\emph{(geometrically) equivalent}
(resp.~\emph{closed orbit equivalent})
if there is a symplectic isomorphism $\varphi:V_1\to V_2$, 
inducing an isomorphism $\tilde\varphi:\mathfrak{sp}(V_1)\to\mathfrak{sp}(V_2)$,
such that $\rho_2(\Lg_2)=\tilde\varphi(\rho_1(\Lg_1))$
(resp.~$\varphi$ maps closed orbits of $G_1$ onto closed orbits 
of $G_2$). The (outer) product of 
$\rho_1$ and $\rho_2$ is the algebra $\Lg_1\oplus\Lg_2$ acting on 
$V_1\oplus V_2$; it is a symplectic representation. A symplectic representation
is called \emph{connected} if it is not equivalent to the 
product of two non-trivial symplectic representations. 
Of course, it suffices to prove 
the above theorems for connected representations.

A symplectic representation $\rho:\Lg\to\mathfrak{sp}(V)$ is called
\emph{saturated} if $\rho[\Lg]$ is self-normalizing in $\mathfrak{sp}(V)$. 
Note that every type~2 representation $U\oplus U^*$ has non-trivial 
endomorphisms, namely, $\Lt^1$ acting by $t\cdot(u,u^*)=(tu,-tu^*)$.
We will also use the following notation from~\cite{Knop1}.
Let $U$ be a representation of a semisimple Lie algebra $\Ls$. 
We denote the type~2 representation of $\Lg=\Ls+\Lt^1$ on $U\oplus U^*$ by $T(U)$. 
Continuing, if $U_1$, $U_2$ are two representations of $\Ls$, then
$T(U_1)\oplus T(U_2)$ is a representation of $\Lg=\Ls+\Lt^2$. 

\begin{rmk}\label{T(U)}
Let $U$ be a symplectic representation of $G$ (so $U\cong U^*$). Then
$(G\times SO_2,U\otimes\C^2)$ is isomorphic to $T(U)=U\oplus U$ via 
$v\otimes e_1+w\otimes e_2\mapsto(v+iw,v-iw)$.
\end{rmk}

Recall that a representation is called \emph{stable} if generic
orbits are closed. A representation of the form
$U\oplus U^*$ is always stable, since it admits the invariant 
orthogonal structure given by 
$\langle u_1+u_1^*,u_2+u_2^*\rangle = u_1^*(u_2)+u_2^*(u_1)$
and one can apply~\cite[Cor.~5.9]{S} or~\cite{Lu2,Lu3}.
A useful necessary and sufficient condition for the stability
of a symplectic representation is that the generic isotropy 
algebra be reductive~\cite[Thm.~2]{Lo}. Recall also that the \emph{rank} of 
a representation $V$ of $G$ is the difference between the dimension 
of $V\cat G$ and that of the subspace of fixed points $V^G$.

%We can now prove:

\begin{prop}\label{commut}
Let $\rho:\Lg\to\mathfrak{sp}(V)$ be a polar symplectic representation. 
%Then:
%\begin{enumerate}
%\item[(a)] The centralizer of $\rho[\Lg]$ in $\mathfrak{sp}(V)$ is commutative. 
%\item[(b)] 
Let $\hat{\Lg}$ be the normalizer of $\rho[\Lg]$ in $\mathfrak{sp}(V)$,
and let $\hat G$ be the corresponding connected subgroup of $Sp(V)$.
Then $(\hat{\Lg},V)$ is saturated and $(\hat G,V)$ is 
closed orbit equivalent to $(G,V)$. 
%\end{enumerate}
It follows %from~(b) 
that $(\tilde{\Lg},V)$ is polar for every 
$\rho[\Lg]\subset\tilde{\Lg}\subset\hat{\Lg}$. 
\end{prop}

\Pf 
%(Compare~\cite[Prop.~2.2]{Knop1}.) 
%Let $V=\bigoplus_i C_i^{n_i}$ be a decomposition into indecomposable
%symplectic representations where the $C_i$ are mutually 
%non-isomorphic. The centralizer of~$\rho[\Lg]$ in~$\mathfrak{sp}(V)$ is
%the product of the centralizers of the $C_i^{n_i}$. There are three cases to consider:
%\begin{enumerate}
%\item[1.] $C_i$ is of type~1. Then the centralizer is $\mathfrak{so}_{n_i}$.
%\item[2a.] $C_i=U\oplus U^*$ is of type~2 with $U\not\cong U^*$. Then the centralizer
%is $\mathfrak{gl}_{n_i}$. 
%\item[2b.] $C_i=U\oplus U^*$ is of type~2 with $U\cong U^*$. Then the centralizer
%is $\mathfrak{sp}_{2n_i}$. 
%\end{enumerate}
%A component $C_i$ of type~$1$ has multiplicity $n_i\leq2$ since $C_i^2$ is stable and 
%we can apply~\cite[Cor.~2.15]{DK} to~$C_i^3$.
%The same corollary yields that $n_i\leq1$ in case 2a since $C_i$ is stable in that case, 
%and that components of type 2b cannot occur since $U$ is stable in that case.
%This proves part~(a).
Since $\Lg$ is reductive, $\hat{\Lg}$ is generated by $\Lg$ and 
its centralizer in $\mathfrak{sp}(V)$.
It follows from Theorem~\ref{polar-implies-coisotropic} 
and~\cite[Lem.~4.1]{Knop1} that the centralizer of $\rho[\Lg]$ in 
$\mathfrak{sp}(V)$ is commutative. Owing to Remark~\ref{T(U)}
and the description of the centralizer in~\cite[Prop.~2.2]{Knop1}, 
we can now write $V=W\oplus U \oplus U^*$, where $W=W_1\oplus\cdots\oplus W_r$,
$U=U_1\oplus\cdots\oplus U_s$, the $W_i$ are indecomposable of type~1,
and either the $U_j\oplus U_j^*$ are indecomposable of type~2
or $U_j$ is of type~1; moreover,  
$\hat{\Lg}=\Lt^s+\Lg=\Lt^s\oplus\Lg'$, 
where $\Lg'$ is the derived algebra of $\Lg$. 

We may assume $U\neq\{0\}$. 
Since $U_j\oplus U_j^*$ is stable, it follows from~\cite[Cor.~2.15]{DK}
that
\[ \mathbb C[V]^G=\mathbb C[W]^G\otimes\mathbb C[U_1\oplus U_1^*]^G\otimes\cdots\otimes\mathbb C[U_s\oplus U_s^*]^G. \]
The final argument in~\cite{DP} (see also~\cite[p.~47]{B-R}) shows that 
$\mathbb C[U_j\oplus U_j^*]^G=\mathbb C[U_j\oplus U_j^*]^{G\cdot T^1}$.
Since $T^s$ acts trivially on $W$, this implies that 
\[ \mathbb C[V]^G=\mathbb C[W]^G\otimes\mathbb C[U\oplus U^*]^{G\cdot T^s}
=\mathbb C[V]^{G\cdot T^s} \]
and the result follows. \EPf

\subsection{Knop reduction}\label{knop-red} 

Fix a Cartan subalgebra $\Lh$ of $\Lg$ and a system 
of positive roots $\Delta^+\subset\Delta$. 
For each $\alpha\in\Delta$, the corresponding coroot
is denoted by $\alpha^\vee$. 
The weight system 
of $V$ is denoted by~$\Lambda$. 
A weight $\lambda\in\Lambda$ 
is called:
\begin{enumerate}
\item[(i)] \emph{extremal} or \emph{highest} if $\alpha\in\Delta$ and 
$\langle\lambda|\alpha^\vee\rangle>0$ implies $\lambda+\alpha\not\in\Lambda$;
\item[(ii)] \emph{toroidal} if $\langle\lambda|\alpha^\vee\rangle=0$
for all $\alpha\in\Delta$;
\item[(iii)] \emph{singular} if it is extremal and $2\lambda\in\Delta$ 
and the multiplicity of $\lambda$ is one. 
\end{enumerate}

A submodule $U$ of~$V$ generated by a highest weight vector
is called \emph{singular} if $U$ is an anisotropic 
subspace of $V$ and $G\to Sp(U)$ is surjective. 
Note that if $\lambda$ is an extremal weight of $V$ 
and $2\lambda\in\Delta$, then 
we can always find a highest weight vector for $\lambda$ that
generates a singular submodule of~$V$; however, in case the 
multiplicity of~$\lambda$ is bigger than one, one can also 
find a highest weight vector that generates an isotropic, hence
non-singular submodule~\cite[Remarks, p.~228]{Knop2}. 

A symplectic representation is called \emph{terminal} if all 
of its highest weights are either toroidal or singular. 
Equivalently, a symplectic representation is terminal 
if every highest weight vector generates either a one-dimensional 
module or a singular submodule. Such a 
representation $(G,V)$ decomposes as
\begin{equation}\label{eq:terminal}
 V=V_0\oplus V_1\oplus \cdots \oplus V_s, \qquad
G= G_0\times Sp(V_1) \times \cdots \times Sp(V_s) 
\end{equation}
where $V_0=\oplus_{i=1}^m(\mathbb C_{\lambda_i}\oplus\mathbb C_{-\lambda_i})$
is a direct sum of $1$-dimensional $G_0$-modules~\cite[Proposition~4.1]{Knop2}.
A terminal symplectic representation is coisotropic if and only if 
the set of weights $\{\lambda_1,\ldots,\lambda_m\}$ 
is linearly independent~\cite[Thm.~3.1]{Knop1}. 

Knop reduction is a finite algorithm which, for a given symplectic 
representation~$(G,V)$, outputs a 
terminal symplectic representation. Indeed if~$(G,V)$ is not itself
terminal, a step of the algorithm is performed by choosing 
extremal weight $\lambda\in\Lambda$ which 
is neither toroidal nor singular and putting 
$P=\{\alpha\in\Delta|\langle\lambda|\alpha^\vee\rangle>0\}$ and $Q=\lambda-P$ 
as multisets (i.e.~sets with multiplicities), and
\[ \Delta'=\Delta\setminus(P\cup-P),\qquad
\Lambda'=\Lambda\setminus(Q\cup-Q). \]
The choice of $\lambda$ ensures that $\Delta'$ is the root system of a
reductive Lie algebra $\Ll$ (namely, a Levi subalgebra of 
the stabilizer of the line through a highest weight vector of $\lambda$),
and $\Lambda'$ is a 
weight system of a symplectic representation $S$ of $\Ll$. 
The main point is that 
$(\Lg,V)$ is coisotropic if and only if so 
is~$(\Ll,S)$~\cite[Thm.~8.4 and Prop.~9.1]{Knop2}. 
This algorithm was used in~\cite{Knop1}
to classify coisotropic symplectic representations. 

\subsection{Relation to polar representations}

\begin{prop}\label{knop-red-polar}
If $(\Lg,V)$ is a non-terminal polar symplectic representation, then 
any Knop reduction $(\Ll,S)$ is also a polar symplectic representation. 
Moreover, any Cartan subspace of $(\Ll,S)$ is a Cartan subspace of
$(\Lg,V)$. 
\end{prop}

\Pf Let $\lambda$ be a highest weight which is neither toroidal  
nor singular. Take a corresponding highest weight vector $v_\lambda$ of unit 
length that generates a non-singular submodule. Consider:
\[ \begin{array}{rl}
v_{-\lambda}=-\epsilon(v_\lambda):&\mbox{lowest weight vector, so
that $\omega(v_\lambda,v_{-\lambda})=1$}\\
\Lp:&\mbox{stabilizer of $\mathbb C\,v_\lambda$ (parabolic subalgebra of $\Lg$)}\\
\Lp_u:&\mbox{unipotent radical of $\Lp$}\\
\Ll:&\mbox{Levi subalgebra of $\Lp$, so that $\Lp=\Ll+\Lp_u$} \\
\Lp^-=\Ll+\Lp_u^-:&\mbox{opposite parabolic subalgebra}
\end{array} \]

Let $v=v_\lambda+v_{-\lambda}$. Then $v$ is a semisimple point~\cite[Proposition~1.2]{DK},
and we may assume it is of minimal length~\cite[pp.~508-509]{DK}.
Now~(\cite[Eq.~(3.3)]{Knop2} or~\cite[Thm.~3.2]{Knop1})
%\begin{equation}\label{s}
\[ S=(\Lp_uv_{-\lambda})^{\omega}\cap(\Lp_u^-v_\lambda)^{\omega}=(\Lg\cdot v)^\omega\oplus \mathbb C\,v. \]
%\end{equation}
Let $c\subset V$ be a Cartan subspace containing~$v$. 
Since 
\[ \omega(c,\Lp_u^-v_\lambda)=h(c,\epsilon\Lp_u^-v_\lambda)=h(c,\Lp_uv_{-\lambda})
=h(c,\Lp_uv)=0, \]
(cf.~\cite[Rem.~1.4]{DK})
and similarly $\omega(c,\Lp_uv_{-\lambda})=0$, we see that $c\subset S$ and
$c=c\cap(\Lg\cdot v)^\omega\oplus\mathbb C\,v$.
We claim that for any regular $x\in c$,
\[ \Ll\cdot x=\Lg\cdot x\cap S. \]
Indeed the direct inclusion is obvious. Moreover, since the
$h$-orthocomplements $N_x$ and $N_v$ to $\Lg\cdot x$ and $\Lg\cdot v$,
resp., satisfy $N_x\subset N_v\subset S$, we have:
\begin{eqnarray*}
\dim(\Lg\cdot x\cap S) & = & \dim S+\dim\Lg\cdot x-\dim V \\
  &=& \dim(\Lg\cdot v)^\omega+1+\dim\Lg\cdot x -\dim V \\
  & = & \dim \Lg\cdot x-\dim\Lg\cdot v + 1\\
  &=& \dim\Lg_v-\dim\Lg_x+ 1\\
 &=& \dim\Ll_v-\dim\Ll_x + 1 \qquad\mbox{(since $\Lg_x\subset\Lg_v\subset\Ll$)}\\
 &=& \dim\Ll\cdot x-\dim\Ll\cdot v+1 \\
 &=&\dim\Ll\cdot x, 
\end{eqnarray*}
which checks the claim. Using $\Lg\cdot c=\Lg\cdot x$, we now
deduce that $\Ll\cdot c=\Ll\cdot x$. 
The proof is finished by noting that 
$\dim c=\dim V\cat G=\dim S\cat L$~\cite[Thm.~8.4]{Knop2}. \EPf

\section{Polar symplectic representations are coisotropic}

In this section we prove Theorem~\ref{polar-implies-coisotropic}.
The proof is along lines suggested by the referee. 
Let $(G,V)$ be a polar symplectic representation. 
We may assume there are no trivial components.

\begin{lem}\label{cartan-subspace-is-isotropic}
Every Cartan subspace is isotropic. 
\end{lem}

\Pf Let $c\subset V$ be a Cartan subspace. The restriction 
$\omega|_{c\times c}$ is $W(c)$-invariant, where
$W(c)=N_G(c)/Z_G(c)$ is
the Weyl group of $(G,V)$ with respect to~$c$, and 
$W(c)$ is generated by unitary reflections~\cite[Lem.~2.7 and Th.~2.10]{DK}. 
For $w\in W(c)$, a vector $u$ in the fixed
point set $c^w$ of~$w$, 
and a $w$-eigenvector $v\in c$ transversal to $c^w$, we have
$\omega(u,v)=\omega(w\cdot u,w\cdot v)=e^{\frac{2\pi i}q}\omega(u,v)$ for some
positive integer $q\neq1$, thus $\omega(c^w,v)=0$. We deduce that 
$v\in\ker\omega|_{c\times c}$. Since a basis of $c$ can be constructed 
which consists of such eigenvectors of reflections
(otherwise $c$ has a non-zero $W(c)$-fixed subspace which implies
that $V$ has a non-zero $G$-fixed subspace~\cite[Lem.~2.11 and
Prop.~2.13]{DK}), 
this shows that the restriction of~$\omega$ to~$c$ is null. \EPf

\medskip

\textit{Proof of Theorem~\ref{polar-implies-coisotropic}.}
Let $v\in V$ be a 
regular element. We may assume $v$ is of minimal length.
Let~$c\subset V$ be a Cartan subspace
containing~$v$. 
Consider first the case in which $V$ is stable. 
Since $v$ is of minimal length,
we have the $h$-orthogonal direct sum $V=c\oplus\Lg\cdot v$. 
Due to~(\ref{omega-h}), $(\Lg\cdot v)^\omega=\epsilon(c)$, and 
by Lemma~\ref{cartan-subspace-is-isotropic}, $\epsilon(c)\subset\Lg\cdot v$,
so $\Lg\cdot v$ is coisotropic. 

In the general case, we use~\cite[Cor.~2.5]{DK} to write $V=c\oplus\Lg\cdot v\oplus U$
where $U$ is a $G_v$-invariant subspace and $U\cat G_v=\{0\}$. 
As above, $\Lg\cdot v\oplus U$ is coisotropic and we need to show that 
it is equal to $T_{(v,u)}G(v,u)$ where $u\in U$ is a generic point. It suffices 
to show that $T_uG_v(u)=U$, or that $G_v$ has an open orbit in $U$. 
Since the action on $\Lg\cdot v=\Lg/\Lg_v$ is orthogonal ($\Lg$ is reductive)
and that on $c$ is trivial, we deduce that $(G_v,U)$ is self-dual. Since $G_v$ has
no nonzero closed orbits in $U$, no component of $(G_v,U)$ is orthogonal and 
$U$ is a sum of pairwise inequivalent indecomposable symplectic representations 
of type~$1$. In particular, the center of $G_v$ acts trivially on $U$. 
Now Theorem~3.3 on p.~165 and the Corollary on p.~156 in~\cite{PV} say that 
the field of fractions of $\mathbb C[U]^{G_v}$ consists of constants only, and $G_v$ 
has an open orbit in $U$. \EPf

\section{The classification}

In this section, we prove Theorem~\ref{classif}.
Thanks to Theorem~\ref{polar-implies-coisotropic}, we will extract
the list of saturated polar symplectic representations from the lists of 
saturated coisotropic representations given, up to geometric 
equivalence, by~\cite[Thms.~2.4, 2.5 and~2.6]{Knop1}. 

%More precisely,
%\cite{Knop1} describes the saturated polar representations as follows:
%Tables~1 and~2 contain the indecomposable coisotropic 
%representations of types~1 and~2, resp.; Tables~11, 12 and 22 contain
%the decomposable ones without $\mathfrak{sl}_2$-links; 
%there is an additional Table~S from which representations 
%with $\mathfrak{sl}_2$-links are constructed
%according to the recipe in Thm.~2.6, \emph{ibid}.

Suppose $V$ is a saturated indecomposable
polar symplectic representation of $\Lg$. If it is 
of type~1, then 
it is listed in~\cite[Table~1]{Knop1}. 
Representations in this table with $\dim V\cat G\leq1$ are trivially polar,
so we run through the other cases.
Some representations with $\dim V\cat G=2$ are already discussed
in~\cite[p.~512 and~523]{DK}.  
We finish this case by referring to~\cite[Tabelle, p.199 and p.201]{L}, where 
irreducible polar representations of connected semisimple Lie 
groups are classified (see also discussion in~\cite[p.~208]{L}).
We obtain our Table~A.

Suppose now $V$ is of type~2. Then it is listed in~\cite[Table~2]{Knop1}. 
In this case $V=U\oplus U^*$ and $(\Lg,V)$ is the complexification of $(\Lk,U)$, 
where $\Lk$ is a maximal compact subalgebra of $\Lg$ and 
$U$ is a real irreducible polar 
representation with an invariant complex structure. Therefore we can refer to the 
classification of irreducible polar representations of compact connected
Lie groups~\cite{D,EH1}. We obtain our Table~B. 

We will complete the proof of the theorem by showing that 
every saturated decomposable polar symplectic representation is a
product,
namely, connected saturated decomposable polar symplectic
representations do not exist. An $\mathfrak{sl}_2$-\emph{link}
is an $\mathfrak{sl}_2$-factor of $\Lg$ which acts effectively on at least 
two indecomposable components of $V$. All connected saturated decomposable
coisotropic representations without 
$\mathfrak{sl}_2$-links are listed 
in~\cite[Tables~11, 12 and~22]{Knop1}, and we will see shortly that 
none of these
is polar. Indeed due to~\cite[Prop.~2.14]{DK}, 
we need only examine the representations in tables~11, 12 and~22 
whose irreducible components are all polar; moreover,
if one of the summands is stable, the \emph{rank condition} 
says that the rank of the sum equals the sum of the ranks of the 
summands. The only unstable representations in Table~A are
$(Sp_{2m},\C^{2m})$ for all $m\geq1$ 
and $(SO_p\times Sp_{2m},\C^p\otimes\C^{2m})$
where $3\leq p < 2m$ and $p$ is odd, while all representations
in Table~B are stable. Now all representations 
in Tables~12 and~22 have both components polar and at least 
one component stable, and we check that the rank condition 
is violated by all of them.
The same argument applies to the representations 
of Table~11, but $\langle11.13\rangle$ which has a non-polar
component and therefore is not polar, and the two sub-cases not having 
stable components of~$\langle11.11\rangle$
and~ $\langle11.14\rangle$, which are discussed in 
Lemmata~\ref{special-case-1} and~\ref{special-case-2}. 

We borrow more notation from~\cite{Knop1} (cf.~(2.4), p.~538).
The line under the $\oplus$-sign below means that the algebras 
immediately to the left and to the right are being identified 
and the resulting algebra is acting diagonally.

\begin{lem}\label{special-case-1}
$\mathfrak{so}_p\otimes\mathfrak{sp}_{2m}\underline{\oplus}\mathfrak{sp}_{2m}$ 
is not polar for $3\leq p<2m$ and $p$ odd. 
\end{lem}

\Pf We will use Proposition~\ref{knop-red-polar}. The 
Lie algebra is $\g{so}_p+\g{sp}_{2m}$ and the representation 
space is $V_1\oplus V_2$, where $V_1=\C^p\otimes\C^{2m}$ and
$V_2=\C^{2m}$. 
By performing Knop reduction with respect to a highest weight vector
of $V_1$
%taking a slice representation at $v_1=v_{\lambda_1}+v_{-\lambda_1}$,
%where $\lambda_1$ is the highest weight of $V_1$, 
and proceeding by induction, we may assume $p=3$ and $m\geq2$. 
A further step of Knop's algorithm yields
\[ \C_{\lambda_1}\oplus\C_{-\lambda_1}\oplus\C_{\epsilon'_1}\oplus\C_{-\epsilon'_1}\oplus
\g{sp}_{2m-2}\underline{\oplus}\mathfrak{sp}_{2m-2}\]
%A further slice representation is $\Lt^1 \oplus \mathfrak{sp}_{2m-2}$ acting 
%on $(\C\oplus\C^*)\oplus(\C^{2m-2}\oplus\C^{2m-2})$ plus a trivial 
%(one-dimensional) component. 
where $\lambda_1=2\epsilon_1+\epsilon'_1$. 
This representation is polar~\cite[p.~522]{DK} 
with Cartan subspace $c=c_0\oplus c_1\oplus c_2$, where 
$c_0$, $c_1$ and $c_2$ are one-dimensional Cartan subspaces
for $\C_{\lambda_1}\oplus\C_{-\lambda_1}$, 
$\C_{\epsilon'_1}\oplus\C_{-\epsilon'_1}$ and
$\g{sp}_{2m-2}\underline{\oplus}\mathfrak{sp}_{2m-2}$, respectively. 
If the given representation were polar, then it would have~$c$ 
as a Cartan subspace.
Since $\C_{\epsilon'_1}\oplus\C_{-\epsilon'_1}\subset V_2$ and
$V_2$ contains no non-zero semisimple points, this is not possible. \EPf 

\begin{lem}\label{special-case-2}
$\mathfrak{sp}_{2m}\otimes\mathfrak{so}_5\underline{\oplus}\mathfrak{sp}_4$ 
is not polar for $m\geq3$.
\end{lem}

\Pf This representation has rank~$4$. Knop reduction with respect to a 
highest weight $\lambda_1$ of the first summand yields
\[ \C_{\lambda_1}\oplus\C_{-\lambda_1}\oplus\mathfrak{sp}_{2m-2}\otimes S^2\mathfrak{sl}_2\underline{\oplus} T(\mathfrak{sl}_2).\]  
Consider the last two summands, namely, 
$\langle S.13\rangle+\langle S.10\rangle$
in~\cite[Table~S]{Knop1}. This is not polar, since its rank is $3$, 
$T(\mathfrak{sl}_2)$ is stable of rank~$1$, and 
$\mathfrak{sp}_{2m-2}\otimes S^2\mathfrak{sl}_2$ has rank $1$, so
the rank condition is not satisfied.\mbox{ }\hfill\mbox{ } \EPf

\medskip

We finish the proof by considering a connected saturated decomposable
coisotropic representation $V$ of $\Lg$   with
$\mathfrak{sl}_2$-links and showing that it cannot be polar.
According to~\cite[Thm.~2.6]{Knop1}, $V$ is
obtained by taking any collection of representations
from Table~S (\emph{ibid}) and identifying any number of disjoint pairs
of underlined $\mathfrak{sl}_2$'s, except that not allowed
is the identification of the two $\g{sl}_2$'s of $\langle S.1\rangle$
and the combination of $\langle S.9\rangle$ with itself. 
Again we need only consider 
entries in Table~S which are polar; for convenience, we list them in 
Table~S'. Note that the only unstable representations therein
are $\langle S.9\rangle$ and $\langle S.13\rangle$ with $m\geq2$. 
An easily checked, common feature of all representations in the table is that
replacing an underlined $\mathfrak{sl}_2$ by $\mathfrak{so}_2$ 
increases the rank by one; we will use this fact below.

\[\begin{array}{|lll|}
\hline
&(\Lg,V) &\\
\hline
\langle S.1\rangle& \underline{\g{sl}}_2\otimes\g{sp}_{2m}\otimes
\underline{\g{sl}}_2 & m\geq1\\
\langle S.3\rangle& \g{so}_n\otimes
\underline{\g{sl}}_2 & n\geq3\\
\langle S.5\rangle& \g{spin}_9\otimes
\underline{\g{sl}}_2 &\\
\langle S.7\rangle& \g{spin}_7\otimes
\underline{\g{sl}}_2 &\\
\langle S.9\rangle& 
\underline{\g{sl}}_2 &\\
\langle S.10\rangle& 
T(\underline{\g{sl}}_2) &\\
\langle S.11\rangle& 
T(\g{sl}_m\otimes\underline{\g{sl}}_2) &m\geq2\\
\langle S.13\rangle& 
\g{sp}_{2m}\otimes S^2\underline{\g{sl}}_2 & m\geq1\\
\langle S.16\rangle& \g{g}_2\otimes
\underline{\g{sl}}_2 &\\
\hline
\end{array} \]
\begin{center}
\textsc{Table S': Polar representations in~\cite[Table~S]{Knop1}}
\end{center}

In view of~\cite[Prop.~2.14]{DK},
we may assume that $V$ has two indecomposable components and both are polar. 
Now $\Lg=\Lg_1+\Ls+\Lg_2$,
$V=V_1\oplus V_2$, and $V_i$ is an indecomposable symplectic
representation of $\Lg_i+\Ls$ given by Table~S', $i=1$, $2$, where 
$\Ls$ is an $\g{sl}_2$-link. Assume at least one of $V_i$ is stable,
say $V_1$. Take the semisimple point $v=v_{\lambda_1}+v_{-\lambda_1}\in V_1$ where
$v_{\lambda_1}$ is a highest weight vector of $V_1$. 
Then the isotropy algebra
$\Lg_{v_1}=(\Lg_1+\Ls)_{v_1}+\Lg_2=(\Lg_1)_{v}+\Lt^1+\Lg_2$.
If $V$ were $G$-polar then, due to~\cite[Cor.~2.15]{DK},  
the set of closed orbits of $T^1\cdot G_2$ in $V_2$ would have to coincide
with the set of closed orbits of $SL_2\cdot G_2$, but it 
follows from the fact above that this is not the case.
The next lemma deals with the remaining two
cases with no stable components and finishes the proof of Theorem~\ref{classif}. 

\begin{lem}
The combination of $\langle S.13\rangle$ with itself or $\langle S.9\rangle$ 
is not polar. 
\end{lem}

\Pf Write $V=V_1\oplus V_2$ where $V_1$ is $\langle S.13\rangle$
and consider Knop reduction $(\Ll,S)$ with respect to the highest weight
$\lambda_1$ of $V_1$. If $V_2$ is $\langle S.13\rangle$, then
$S$ contains as a summand
$T(\g{sp}_{2m})\underline{\oplus}\g{sp}_{2m}$
which has already been shown not to be polar.
Hence $V$ is not polar.

On the other hand, if $V_2$
is $\langle S.9\rangle$, then $S$ is the polar representation
$(\C_{\lambda_1}\oplus\C_{-\lambda_1}\oplus U)\oplus(\C_{\epsilon_1}\oplus\C_{-\epsilon_1})$, where $U$ is a subspace of $V_1$, $\pm2\epsilon_1$ are the roots of
$\g{sl}_2$, and $\C_{\epsilon_1}\oplus\C_{-\epsilon_1}$ equals $V_2$.
By~\cite[Prop.~2.14]{DK}, a Cartan subspace
of this representation is of the form $c=c_1\oplus c_2$, where 
$c_1\subset\C_{\lambda_1}\oplus\C_{-\lambda_1}\oplus U$ and $c_2$ is the diagonal 
subspace of $\C_{\epsilon_1}\oplus\C_{-\epsilon_1}$.
If~$V$ were polar, 
Proposition~\ref{knop-red-polar} says that 
$c$ would be a Cartan subspace
of $V$. However $V_2$ does not contain
non-zero $G$-semisimple elements of $V$, so this is not
possible.\EPf

%\begin{lem}\label{13-13}
%The combination $\langle S.13\rangle+\langle S.13\rangle$ is not polar. 
%\end{lem}
%
%\Pf This representation is $\g{sp}_{2m}\otimes S^2\g{sl}_2 \underline{\oplus}
%S^2\g{sl}_2\otimes\g{sp}_{2m}$. A slice representation is given by
%$T(\g{sp}_{2m})\underline{\oplus}\g{sp}_{2m}\oplus\g{sp}_{2m}$. 
%However $T(\g{sp}_{2m})\underline{\oplus}\g{sp}_{2m}$ is not polar 
%since it has rank~$2$,  $T(\g{sp}_{2m})$ is stable of rank~$1$, 
%$\g{sp}_{2m}$ has rank zero, and we can apply Proposition~\ref{decomp}(b). 
%The result now follows from Proposition~\ref{slice}. \EPf

\section{Symplectic symmetric spaces}

A \emph{symplectic symmetric space} is a symmetric space which is endowed
with a symplectic structure invariant by the symmetries. We refer
to~\cite{Bie0,Bie1} for the basic theory of such spaces. 
Our interest in them is that the (complexified) isotropy representations of 
symplectic symmetric spaces provide examples of symplectic 
\emph{$\theta$-groups} (namely, adjoint groups of graded 
Lie algebras) thus, polar symplectic 
representations~\cite[\S8.5, 8.6]{PV}.  
Indeed simply-connected
symplectic symmetric spaces are parametrized by symplectic involutive 
Lie algebras. A \emph{symplectic involutive Lie algebra} is a
triple~$(\Lg,\sigma,\omega)$ where $\Lg$ is a real Lie algebra,
$\sigma$ is an involution of $\Lg$, with respect to which there is
an eigenspace decomposition $\Lg=\Lh+\Lq$, and $\omega$ is an 
$\mathrm{ad}_{\mathfrak h}$-invariant non-degenerate $2$-form on $\Lq$. 

An indecomposable (i.e.~non-isomorphic to a product of symplectic involutive
Lie algebras) non-flat (i.e.~satisfying $[\Lq,\Lq]\neq0$) 
reductive symplectic involutive Lie algebra $(\Lg,\sigma,\omega)$
is simple~\cite[Prop.~3.5.4]{Bie0}. The symplectic structures $\omega$ on a 
simple involutive Lie algebra $(\Lg,\sigma)$ are parametrized
by the non-zero elements in the center $Z(\Lh)$ of $\Lh$~\cite[Th.~2.1]{Bie1}.
Moreover, if $\Lg$ is a complex Lie algebra viewed as real, then
$\sigma$ is a complex automorphism, $\omega$ is complex bilinear,
$\dim_{\mathbb C} Z(\Lh)=1$ and $(\Lh,\Lq)$ is a $\theta$-group; 
otherwise $\Lg$ is absolutely simple, $\dim_{\mathbb R} Z(\Lh)=1$
and the complexification $(\Lh^{\mathbb C},\Lq^{\mathbb C})$ is a 
$\theta$-group~\cite[Prop.~2.2 and Thm.~2.2]{Bie1}. In any case,
the (indecomposable) polar symplectic representations thus obtained are exacly the 
complexified isotropy representations of irreducible Hermitian Riemannian 
symmetric spaces~\cite[\S~10]{Bie1}, hence of type~2 and listed in Table~B. 
On the other hand, every representation
in Table~B is closed orbit equivalent to the complexification
of the isotropy representation of an irreducible Hermitian Riemannian 
symmetric space~\cite{D,EH1}. 

A complex semisimple symplectic involutive Lie algebra is the 
product of complex simple symplectic involutive Lie 
algebras~\cite[Prop.~3]{BCG}, each of which with an indecomposable 
(polar) symplectic representation (of type~2) as isotropy representation. 
Now an arbitrary polar symplectic representation can be assumed saturated,
up to closed orbit equivalence (Proposition~\ref{commut}), and then it is
the product of indecomposable polar symplectic representations 
(Theorem~\ref{classif}); finally, it is closed orbit equivalent to
the isotropy representation of a complex semisimple symplectic symmetric 
space if and only if each of its factors is. 
This completes the proof of Theorem~\ref{symmetric-spaces}. 

\section{The moment map}\label{sec:moment}

Recall the idea of a moment map. In our context, an action of an algebraic
group $G$ on a symplectic variety $(X,\omega)$ is called \emph{Hamiltonian}
if there exists a \emph{moment map}, that is, an equivariant map
$\mu:X\to\Lg^*$ (where $\Lg^*$ is regarded with the coadjoint 
representation) such that $\omega(\xi x,v) = \langle d\mu_x(v)|\xi\rangle$
for all $\xi\in\Lg$, $v\in T_xX$, $x\in X$ (compare~\cite[section 2]{Knop2}).  
In our particular case of interest $X=V$ is a symplectic representation of 
$G$, there is a canonical moment map given by 
\[ \mu:V\to\Lg^*,\qquad \langle\mu(v)|\xi\rangle=\frac12\omega(\xi v,v). \]

Assume now $(G,V)$ is a saturated polar symplectic representation. 
Apply Knop reduction to get a terminal representation~(\ref{eq:terminal})
with set of 
weights $\{\lambda_1,\ldots,\lambda_r\}$. 
Let $v_{\lambda_j}$ be an $h$-unit $\lambda_j$-weight vector, and 
$v_{-\lambda_j}=-\epsilon(v_{\lambda_j})$
so that $\omega(v_{\lambda_j},v_{-\lambda_j})=1$. 
We easily see from Proposition~\ref{knop-red-polar} that 
$c=\langle v_{\lambda_1}+v_{-\lambda_1},\ldots,  v_{\lambda_r}+v_{-\lambda_r}\rangle$
is a Cartan subspace of $(G,V)$.

%\Pf Note that 
%$\dim c = \dim V\cat G$ in view of Theorem~\ref{polar-implies-coisotropic}
%and~\cite[Th.~3.1]{Knop1}.

%We proceed to the proof by induction on $m(\Lg,V)$, the initial 
%case $m(\Lg,V)=0$ being clear (cf.~Remark~\ref{rmk:terminal}). 
%Let $\lambda=\lambda_1$ be a non-toroidal, non-singular extremal 
%weight and perform one step of Knop reduction to get $(\Ll,S)$. 
%Since $(\Ll,S)$ and $(\Lg,V)$ have the same Knop reduction to a 
%terminal representation and $m(\Ll,S)<m(\Lg,V)$, $c$ is a 
%Cartan subspace for $(\Ll,S)$. Consider now 
%$v=v_\lambda+v_{-\lambda}\in S$. Then $\Ll\cdot v= \mathbb C(v_\lambda-
%v_{-\lambda})$, so the $h$-orthocomplement of $\Ll\cdot v$ in $S$ 
%is $N_v^S=S_\circ\oplus\mathbb C\, v$. Note that $c\perp_h v_\lambda-v_{-\lambda}$,
%so $c\subset N_v^S$.
%By the proof of~\cite[Th.~2.4]{DK}, $c$ is a Cartan subspace 
%for $(L_v,N_v^S)$.  Since $\Ll_v=\Lg_v$ (\ref{l-gv}),
%$N_v^S\subset N_v$ and $\dim c=\dim N_v^S\cat G_v=\dim S\cat L=\dim V\cat G=\dim N_v\cat G_v$, it follows that $c$ is a Cartan subspace for $(G_v,N_v)$.
%Another application of~\cite[Th.~2.4]{DK} finally yields that 
%$c$ is a Cartan subspace for $(G,V)$. 
%\EPf

%\medskip

%We say that a set of weights of a representation is 
%\emph{strongly orthogonal} if neither the sum nor the difference of weights 
%in the set is a root. 

\begin{prop}\label{moment}
\begin{enumerate}
\item[(a)] The set $\{\lambda_1,\ldots,\lambda_r\}$ is strongly
orthogonal in the sense that $\lambda_i+\alpha\neq\pm\lambda_j$
for all $i\neq j$ and $\alpha\in\Delta$. 
\item[(b)] 
  \begin{equation}\label{mu}
    \mu\left(\sum_{j=1}^r a_j(v_{\lambda_j}+v_{-\lambda_j})\right)
    =\sum_{j=1}^r a_j^2\lambda_j,
    \end{equation}
where $a_j\in\mathbb C$. 
In particular,
\[ \mu(c) = \langle\lambda_1,\ldots,\lambda_r\rangle =:\La^*\subset\Lh^* \]
and $\mu$ maps closed orbits to closed orbits.
\end{enumerate}
\end{prop}

\Pf (a) We may assume $V$ is indecomposable and $\mathrm{rank}\,(V)\geq2$. 
One sees that~$V$ is \emph{weight multiplicity free} (wmf), 
in the sense that all the multiplicities of its $\Lh$-weights are at most 
one. Indeed this follows from~\cite[Electronic version]{Knop3} in case $V$ is of type~$2$.
If $V$ is of type~$1$, then~$V$ is polar irreducible and $G$ is semisimple, 
which implies that $V$ is visible~\cite[p.~194]{L} (without zero weights),
and hence wmf due to~\cite[Lem.~3.4]{Kac2}. 

By polarity, $h(\Lg_\alpha(v_{\lambda_j}+v_{-\lambda_j}),v_{\lambda_k}+v_{-\lambda_k})\subset
h(\Lg_\alpha\cdot c,c)=0$ for all $\alpha\in\Delta$ and $j$, $k=1,\ldots,r$.  
Taking components yields
$h(\Lg_\alpha\cdot v_{\pm\lambda_j},v_{\pm\lambda_k})=0$. Since our representation 
is wmf, this implies that $\{\lambda_1,\ldots,\lambda_r\}$ is 
strongly orthogonal. 

(b) $\omega(\Lg_\alpha(v_{\lambda_j}+v_{-\lambda_j}),v_{\lambda_k}+v_{-\lambda_k})=0$ for all $\alpha\in\Delta$, 
by strong orthogonality of $\lambda_j$, $\lambda_k$ in case 
$j\neq k$, and by non-singularity of $v_{\lambda_j}$ in case $j=k$. 
This already shows $\mu(c)\subset\Lh^*$. 
To finish, let $\xi\in\Lh$ and compute
\begin{eqnarray*}
\omega(\xi(v_{\lambda_j}+v_{-\lambda_j}),v_{\lambda_k}+v_{-\lambda_k}) 
&=&\langle\lambda_j|\xi\rangle\,\omega(v_{\lambda_j}-v_{-\lambda_j},v_{\lambda_k}+v_{-\lambda_k})\\
&=&\left\{\begin{array}{ll}
                        0, \mbox{if $j\neq k$},\\
                        2\langle\lambda_j|\xi\rangle, \mbox{if $j=k$}. 
   \end{array}\right. 
\end{eqnarray*}
The desired formula follows. \EPf

\medskip

Since the moment map is equivariant, there is an
induced \emph{invariant moment map}:
\begin{diagram}
V & \rTo^{\mu} & \Lg^* \\
\dTo && \dTo \\
V\cat G & \rDotsto_{\psi:=\mu\cat G} & \Lg^*\cat G\\
\end{diagram}
By Chevalley's theorem, $\Lg^*\cat G\cong \Lh^*/W_G$.
Similarly, by polarity $V\cat G\cong c/W(c)$, where $W(c)=N_G(c)/Z_G(c)$ is
the Weyl group of $(G,V)$ with respect to~$c$; in addition,
since $G$ is connected, $\mathbb C[c]^{W(c)}$ is a polynomial algebra
~\cite[Th.~2.9 and~2.10]{DK}. Hence $V\cat G\cong\C^{\dim c}$.
Let $\pi:\Lh^*\to\Lh^*/W_G$ denote the projection, where 
$W_G$ denotes the Weyl group of $\Lg$ with respect to~$\Lh$.   
By coisotropicity, there is 
a subgroup $W_V$ of $\Gamma:=N_{W_G}(\La^*)/Z_{W_G}(\La^*)$
acting on $\La^*$ as a group generated by reflections such that 
the image of the (invariant) moment map is $\pi[\La^*]$, 
$\psi$ factors through an isomorphism
$V\cat G\to\La^*/W_V$ and 
the morphism $\La^*/W_V\to\pi[\La^*]$ is finite~\cite{Knop2}.  
In view of~(\ref{mu}), $\mathbb C[c]^{W(c)}$ consists of polynomials
in $a_1^2,\ldots,a_r^2$. We deduce there exists a 
$(\mathbb Z_2)^r$-subgroup $D\subset W(c)$ with generators given, for each $i_0$, 
by a transformation that maps
$v_{\lambda_{i_0}}+v_{-\lambda_{i_0}}$ to its opposite and fixes 
$v_{\lambda_{i}}+v_{-\lambda_{i}}$ for all $i\neq i_0$. 

\begin{cor}
The group $W(c)$ is an extension of $W_V$ by $D$.
\end{cor}

\Pf Since $\mu:V\to\Lg^*$ is $G$-equivariant and maps $c$ onto~$\La^*$,
it induces a homomorphism $\bar\mu:W(c)\to W_V$ whose kernel is precisely $D$. 
Moreover $\mu$ induces an isomorphism $c/D\cong\La^*$ and $W(c)/D$ acts
on $c/D$ with quotient isomorphic to $\La^*/W_V$. It follows that 
$\bar\mu$ is onto, namely, $W(c)/D\cong W_V$. \EPf

%\begin{prop}[Knop]
%Let $(G,V)$ be a coisotropic representation. 
%Then the following assertions are equivalent:
%\begin{enumerate}
%\item $\psi$ is an isomorphism onto its image.
%\item $\pi[\La^*]$ is a normal affine variety. 
%\item $W_G\cdot\xi\cap\La^*=\Gamma\cdot\xi$ for all $\xi\in\La^*$.
%\item $\pi[\La^*]\cong\mathbb C^{\dim\mathfrak a}$. 
%\item $\alpha$ is surjective.
%\item The homomorhism $\psi^*:\mathbb C[\Lg^*]^G\to\mathbb C[V]^G$ 
%is surjective. 
%\end{enumerate}
%Moreover, they all imply that $W_V=\Gamma$. 
%\end{prop}
%
%\Pf Note that $\pi$ is a closed map, so $\pi[\La^*]$ is a closed 
%subset of the affine variety $\Lh^*/W_G$ and thus affine. 
%Since $(G,V)$ is assumed coisotropic, $V\cat G\cong\La^*/W_V\cong\mathbb C^{\dim\mathfrak a}$. 

\medskip

\textit{Proof of Theorem~\ref{moment-map}.}
It remains only to prove the second assertion.
Note that $\psi^*:\C[\Lg^*]^G\to\C[V]^G$ factors
as the composition
\begin{diagram}
 \C[\Lh^*]^{W_G} & \rTo^\alpha & \C[\La^*]^\Gamma & \rTo^\beta & \C[\La^*]^{W_V}.
\end{diagram}
We may assume that~$V$ is indecomposable. 
According to the last column of Tables~1 and~2 in~\cite{Knop1}, 
$W_V=\Gamma$ so that $\beta$ is the identity map, and $\alpha$ 
is surjective in all cases but $T(\g e_6)$.
We finish the proof by proving directly in this case 
that $\psi$ is injective. 

We need to show that $W_G\cdot \xi\cap\La^* =\Gamma\cdot\xi$ 
for all $\xi\in\La^*$ in case $V=T(\g e_6)$. We have
$\La^*=\langle\varpi_1+\eta,-\varpi_6+\eta,-\varpi_1+\varpi_6+\eta\rangle$
(we use Bourbaki's notation for weights~\cite[Planches I-IX]{Bourbaki} 
and denote by $\eta$ the weight associated to $\Lt^1$). 
Let $\xi_1$, $\xi_2\in\La^*$ such that $w\xi_1=\xi_2$ 
for some $w\in W_G$. The action of $W_G$ fixes $\eta$ and
%preserves $\tilde{\La}^*=\langle\varpi_1$, $\varpi_6\rangle$, 
$w\tilde\xi_1=\tilde\xi_2$, where $\tilde\xi_i$ 
is the $\tilde{\La}^*$-component of $\xi_i$ with respect 
to~$\La^*=\tilde{\La}^*\oplus\C\eta$.
Let $\Delta_0$ be the root subsystem of type $\sf D_4$ spanned by the simple
roots $\alpha_2,\alpha_3,\alpha_4,\alpha_5$. 
The subsystem $\Delta^i$ of $\Delta$ consisting of roots 
orthogonal to~$\tilde\xi_i$ 
contains $\Delta_0$, so $\Delta^2$ is one of: $\Delta_0$,
all of $\Delta$, or a subsystem of type~$\sf D_5$ containing $\Delta_0$.  
In the first case $\Delta^2=\Delta_0=\Delta^1$, so 
$w\in N_{W_G}(\Delta_0)=N_{W_G}(\tilde{\La}^*)$. Further, any two subsystems 
of type~$\sf D_4$ of a root system of type~$\sf D_5$ are $W(\sf D_5)$-conjugate,
so in the other two cases
there is $w'\in W(\Delta^2)$ such that $w'w(\Delta_0)=\Delta_0$.
Now $w'w\in N_{W_G}(\tilde{\La}^*)$ and $w'w\tilde\xi_1=\tilde\xi_2$. \EPf

\begin{rmk}
Since $\pi:\Lh^*\to\Lh^*/W_G$ is a dominant finite morphism
between affine varieties, it is a closed map and thus $\pi[\La^*]$ 
is an affine variety. We have shown that $\psi$ is a bijective
morphism from $V\cat G$ to $\pi[\La^*]$.    
Essentially by Zariski's main theorem~\cite[ch.~8]{Milne},
$\psi$ is an isomorphism onto its image 
if and only if $\pi[\La^*]$ is a normal variety.
In general, $\psi:V\cat G\cong\La^*/\Gamma\to\pi[\La^*]$ is the 
normalization morphism.
\end{rmk}

\begin{ex}
In Theorem~\ref{moment-map}, we cannot drop the assumption that
the representation is saturated. In fact, consider
the polar symplectic representation of type~2 given 
by~$(SL_n,\C^n\oplus\C^{n*})$. Then $\mu:\C^n\oplus\C^{n*}\to \mathfrak{sl}_n^*$
is given by $\mu(u,\alpha)(\xi)=\alpha(\xi\cdot u)$ and 
$\psi:\C\to\C^{n-1}$ is given by $\psi(z)=(\sigma_2,\ldots,\sigma_n)$
where $\sigma_j=-(j-1){{n}\choose{j}}\left(\frac zn\right)^j$.
In particular, if $n=2$ then $\psi(z)=-\frac{z^2}4$ does not separate closed
orbits. Note that if $n=3$ then $\psi(z)=(-\frac{z^2}3,-\frac{2z^3}{27})$ 
is not 
an isomorphism onto its image; however, for all~$n$,
the enlarged saturated polar 
symplectic representation $T(\g{sl}_n)$ has $\C\to\C^n$, 
$z\mapsto(z,0,\ldots,0)$ as moment map, that is, an isomorphism onto its image. 
\end{ex}

%\bibliographystyle{amsalpha}
%\bibliography{ref}

\providecommand{\bysame}{\leavevmode\hbox to3em{\hrulefill}\thinspace}
\providecommand{\MR}{\relax\ifhmode\unskip\space\fi MR }
% \MRhref is called by the amsart/book/proc definition of \MR.
\providecommand{\MRhref}[2]{%
  \href{http://www.ams.org/mathscinet-getitem?mr=#1}{#2}
}
\providecommand{\href}[2]{#2}

\end{document}